\documentclass[oneside,english,reqno]{amsart}
\usepackage[T1]{fontenc}
\usepackage[latin1]{inputenc}
\usepackage{amssymb}

\makeatletter

 \theoremstyle{plain}
 \theoremstyle{plain}    
 \newtheorem{thm}{Theorem} 
 \theoremstyle{remark}
 \newtheorem*{rem*}{Remark}
 \theoremstyle{plain}    
 \newtheorem*{cor*}{Corollary}

\usepackage{amsfonts}

\usepackage{graphicx}
\usepackage{amscd}

\DeclareMathOperator{\real}{Re}
\DeclareMathOperator{\imag}{Im}
\DeclareMathOperator{\ind}{ind}
\DeclareMathOperator{\spec}{spec}

\theoremstyle{plain}

\usepackage{babel}
\makeatother
\begin{document}

\title{On the gaps between zeros of trigonometric polynomials}

\author{Gady Kozma}

\curraddr{Weizmann Institute of Science}

\email{gadyk@wisdom.weizmann.ac.il, gadykozma@hotmail.com}

\thanks{Supported by EEC Research Training Network {}``Classical Analysis,
Operator Theory, Geometry of Banach spaces, their interplay and their
applications'', contract HPRN-CT-00116-2000.}

\author{Ferenc Oravecz}

\curraddr{Graduate School of Information Sciences, Tohoku University, Sendai,
Japan}

\email{oravecz@renyi.hu}


\begin{abstract}
We show that for every finite set $0\notin S\subset\mathbb{Z}^{d}$
with the property $-S=S$, every real trigonometric polynomial $f$
on the $d$ dimensional torus $\mathbb{T}^{d}=\mathbb{R}^{d}/\mathbb{Z}^{d}$
with spectrum in $S$ has a zero in every closed ball of diameter
$D\left(S\right)$, where\[
D\left(S\right)=\sum_{\lambda\in S}\frac{1}{4||\lambda||_{2}}\quad,\]
 and investigate tightness in some special cases. 
\end{abstract}
\maketitle

\section{Introduction and presentation of the main results}

We are interested in the following question: for real trigonometric
polynomials with a given spectrum, what is the maximal possible distance
between two consecutive real zeros? Of course, if we can use enough
frequencies for building our trigonometric polynomial, the size of
the largest gap can be as close to $2\pi$ (or to $1$, if one uses
the transformation used throughout this paper) as we want. The question
gets more interesting when the number of frequencies is relatively
small. 

The first theorem is a very general result, partially answering the
above question. We state it for multivariate trigonometric polynomials
on the $d$-dimensional torus $\mathbb{T}^{d}=\mathbb{R}^{d}/\mathbb{Z}^{d}$.

\begin{thm}
Let $0\not\in S\subset\mathbb{Z}^{d}$ be a finite set s.t. $-S=S$.
Let \[
f(x)=\sum_{\lambda\in S}c(\lambda)\exp2\pi i\langle x,\lambda\rangle\]
 be a real value trigonometric polynomial on $\mathbb{T}^{d}$. Then
$f$ has a zero in any closed ball of diameter $D(S)$ where\begin{align*}
D(S) & :=\sum_{\lambda\in S}\frac{1}{4||\lambda||_{2}}\quad,\\
||\lambda||_{2} & =\sqrt{\langle\lambda,\lambda\rangle}=\sqrt{{\textstyle \sum_{j=1}^{d}\lambda_{j}^{2}}}\quad.\end{align*}

\end{thm}
Note that the requirement $-S=S$ is necessary because we are interested
in real trigonometric polynomials.

How sharp is this theorem? To answer that we have to analyze specific
cases. From now on we concentrate on the one dimensional situation.
For a finite set $0\notin S\subset\mathbb{Z}$ define \[
M\left(S\right):=\sup_{\spec f\subset S}\sup_{\substack{I\subset[0,1]\\
f|_{I}\neq0}
}|I|\quad,\]
 where $\spec f$ is the spectrum of $f$ i.e.~the support of the
the Fourier transform $\hat{f}$. The second supremum is over all
intervals $I$ such that $f$ is never zero on $I$. Here and also
later on $|I|$ stands for the length of $I$. With this notation
theorem 1 is shortened to the formula\[
M(S)\leq D(S)\quad.\]

The easiest case to investigate is when the spectrum is a (discrete)
interval, or rather, two intervals situated symmetrically around $0$.
It turns out that in this case it is possible to get an explicit formula
for the supremum (maximum does not exist) of the gap size: 

\begin{thm}
\label{thm:interval}For $S=[-N-K,-N]\cup[N,N+K]\subset\mathbb{Z}$
with $N$, $K\in\mathbb{N}$,\[
M(S)=\frac{K+1}{2N+K}\quad.\]

\end{thm}
As the proof will show, we can also demonstrate an explicit formula
for the trigonometric polynomial corresponding to the supremum (being
a polynomial with {}``touching zeros'' --- e.g.~non-negative on
an interval of length $\left(K+1\right)/\left(2N+K\right)$). We see
that the general result applied for this specific case is sharp only
asymptotically when $K\ll N$. This leaves many questions. As examples,
we formulate some of them: 

\theoremstyle{remark}

\newtheorem{question}{Question}

\begin{question}

Theorem 2 can be generalized to linear progressions with stepping
$b$ smaller than twice the starting element, i.e.~when $S=\pm\left\{ N+nb:n\in0,1,2,\dots K\right\} $
for some $K,b\in\mathbb{N}$, $b<2N$ (see Theorem 3 on page \pageref{thm:Mlinprog}).
What happens when $b\geq2N$? Can an explicit formula for $M\left(S\right)$
be found in this case too? 

\end{question}

\begin{question}

What is the maximum when the spectrum points are perfect squares,
i.e.~$S=\left\{ \pm n^{2}\right\} _{n=N}^{N+K}$ ? 

\end{question}

Theorem 1 states that in this case the maximum tends to zero as $N\rightarrow\infty$
even if $K\rightarrow\infty$, a result which by itself is intriguing.
However, is the asymptotic truly $O(1/N)$ or is it even lower? 

\begin{question}

What can we say when $S$ is a random set? A typical model is taking
every integer in $1,\dotsc,N$ with probability $\tau$ and independently
and then symmetrizing. For asymptotic results we should assume $\tau$
and $N$ are related somehow, say $\tau=N^{-\delta}$ for some $0<\delta<1$.

\end{question}

\begin{question}

What about nets? A net of order $n$ is a set with $2^{n+1}$ elements,
defined by parameters $\left\{ a_{0},\dots,a_{n}\right\} $ as follows:
$S:=\pm\left\{ a_{0}+\sum_{i=1}^{n}\left\{ 0,a_{i}\right\} \right\} $.
In many questions nets behave like generalized arithmetic sequences.
Is it true in our case as well? 

\end{question}

We should remark that the specific case of $N=1$ in theorem \ref{thm:interval},
together with the extremal polynomial, has been known for a while.
The question was asked in \cite{T65} and answered in \cite{B84},
which actually showed that this polynomial maximizes the \textbf{measure}
of the positive part. The same polynomial reappears attributed to
Arnold in \cite{T97} and again in \cite{GS00} (which also handles
the question of positive coefficients). For our needs, though, the
case $N=1$ is not enough to make a meaningful comparison between
$M$ and $D$ since we are even more interested in the case when they
are both small, and the proofs in these papers do not seem to generalize
nicely. See also \cite{Y96} for inexact multidimensional results.
Finally, \cite{BH84,MDS99} deal with the related question of completely
positive polynomials with some bounds on the zeroth coefficient.

\section{Proof of Theorem 1 and Theorem 2}

\begin{proof}
[Proof of theorem 1]We will use induction on the cardinality of $S$.
Due to the symmetry of $S$, we have $\# S=2k$, $k=0,1,2,\dots$
The case $k=0$ is trivial as then $f\equiv0$.

In general, given $k>0$, assume we have proven the statement for
all sets $S$ with cardinality smaller than $2k$. Take a set $S$
with cardinality $2k$ and a trigonometric polynomial $f\left(x\right)=\sum_{\lambda\in S}c\left(\lambda\right)\exp2\pi i\left\langle x,\lambda\right\rangle $.
Arguing by contradiction, assume that $f$ does not have a zero, say
is strictly positive, on a $d$-dimensional ball of radius $R=\frac{1}{2}D(S)$
centered at $y=(y_{1},\dotsc,y_{d})$, which we denote as usual by
$B(y,R)$. Pick a frequency $\nu=\left(\nu_{1},\dots,\nu_{d}\right)\in S$,
and translate $f$ by $\mu:=\frac{\nu}{2||\nu||_{2}^{2}}$. Since
$f$ is strictly positive on the ball $B\left(y,R\right)$, the translate
is strictly positive on the translated ball $B(y-\mu,R)$. Adding
the translate to $f\left(x\right)$ we get a new function \[
\tilde{f}\left(x\right):=f\left(x\right)+f\left(x+\mu\right)\]
 which is strictly positive on the closed ball\begin{equation}
B\left(y-{\textstyle \frac{1}{2}}\mu,R-{\textstyle \frac{1}{2}}||\mu||_{2}\right)\subseteq B\left(y,R\right)\cap B\left(y-\mu,R\right)\label{eq:intersect_balls}\end{equation}
 whose diameter is $2R-\frac{1}{2||\nu||_{2}}$.

The Fourier expansion of the translated function is \begin{eqnarray*}
f\left(x+\mu\right) & = & \sum_{\lambda\in S}c\left(\lambda\right)\exp2\pi i\left\langle x+\mu,\lambda\right\rangle \\
 & = & \sum_{\lambda\in S}c(\lambda)\exp2\pi i\langle\mu,\lambda\rangle\exp2\pi i\langle x,\lambda\rangle\end{eqnarray*}
 and thus \[
\tilde{f}\left(x\right)=\sum_{\lambda\in S}c(\lambda)(1+\exp2\pi i\langle\mu,\lambda\rangle)\exp2\pi i\left\langle x,\lambda\right\rangle .\]
Since for $\lambda=\pm\nu$ we get $\exp2\pi i\langle\mu,\lambda\rangle=e^{\pm\pi i}=-1$,
the spectrum $S_{\tilde{f}}:=\spec(\tilde{f)}$ of $\tilde{f}$ is
contained in $S$ but definitely does not contain $\pm\nu$. Thus
$D(S)-\frac{1}{2||\nu||_{2}}\geq D(S_{\tilde{f}})$, and so we find
that we have a new set $S_{\tilde{f}}$ of cardinality smaller than
$2k$ and a trigonometric polynomial $\tilde{f}$ with spectrum $S_{\tilde{f}}$
which is strictly positive on the closed ball $B(y-\mu,R(S_{\tilde{f}}))$
of radius $R(S_{\tilde{f}})=\frac{1}{2}D(S_{\tilde{f}})$ centered
at $y-\mu$. This contradicts our inductive hypothesis.
\end{proof}
\begin{rem*}
A similar proof can be applied to shapes different from balls. For
example, for cubes we get the following version of theorem 1: If $\spec f=S$
then $f$ has a zero in every closed cube of side length $L(S)$ where
\[
L(S):=\sum_{\lambda\in S}\frac{||\lambda||_{\infty}}{4||\lambda||_{2}^{2}}\quad.\]
The only proof element that needs modification is (\ref{eq:intersect_balls}),
as for cubes we have (denoting by $C(y,L)$ the cube of side length
$L$ based at $y$),\[
C(y,L-||\mu||_{\infty})\subset C(y,L)\cap C(y-\mu,L)\]
and $||\mu||_{\infty}=||\nu||_{\infty}/||\nu||_{2}^{2}$.
\end{rem*}

\subsection*{Proof of Theorem 2.}

Theorem 2 is a specific case of the following 

\begin{thm}
\label{thm:Mlinprog}For $S=\{-N-Kb,\dots,-N-b,-N\}\cup\{ N,N+b,\dots,N+Kb\}$
with $N$, $K$, $b\in\mathbb{N}$ and $b<2N$ we have \[
M(S)=\frac{K+1}{bK+2N}\quad.\]

\end{thm}
\begin{proof}
First we show \begin{equation}
M\left(S\right)\leq\frac{K+1}{bK+2N}\label{eq:MltM}\end{equation}
 by proving that for every trigonometric polynomial $f\left(t\right)=\sum_{\lambda\in S}c\left(\lambda\right)e^{2\pi i\lambda t}$,
if $f>0$ on an interval $I$ then $|I|\leq\frac{K+1}{bK+2N}.$ We
may assume w.l.o.g.~that $I=\left[0,a\right]$. One can write $f\left(x\right)=\real F\left(t\right)$
with $F\left(t\right):=e^{2\pi iNt}Q\left(e^{2\pi it}\right),$ where
$Q$ is a polynomial of degree $Kb$, that is\begin{equation}
Q\left(z\right)=c\prod_{j=1}^{Kb}\left(z-\xi_{j}\right)\quad.\label{eq:Q_is_prod}\end{equation}
 The condition $f>0$ on $I$ is equivalent to $\arg F\in\left(-\pi/2,\pi/2\right)$
on $I$. Since $e^{2\pi iNt}$ does a rotation of $2\pi Na$ in the
positive direction as $t$ goes from $0$ to $a$, $f>0$ on $I$
requires that $Q\left(e^{2\pi it}\right)$ does a rotation of at least
$2\pi Na-\pi$ in the negative direction on $I$. In other words,
if $\widetilde{\arg}Q$ is any continuous version of $\arg$ on the
arc $\left\{ e^{2\pi it}:t\in I\right\} $ ($\widetilde{\arg}$ exists
because $\left.f\right|_{I}>0$ gives $Q\neq0$ on the arc, so for
example take $\widetilde{\arg}\, Q\left(e^{2\pi it}\right):=\imag\int_{0}^{t}\frac{Q\prime}{Q}$)
then a necessary condition is \begin{equation}
\ind Q:=\widetilde{\arg}\, Q\left(e^{2\pi ia}\right)-\widetilde{\arg}\, Q\left(1\right)\leq\pi-2\pi Na\quad.\label{eq:indQsmall}\end{equation}
 Clearly the index $\ind Q$ of $Q$ does not depend on the choice
of $\widetilde{\arg}$ (all versions of $\widetilde{\arg}$ differ
by a constant $2\pi k$). As we want an upper bound for $a=|I|$,
we want to maximize $-\ind Q$. Clearly\[
\ind Q=\sum_{j=1}^{Kb}\ind\left(z-\xi_{j}\right)\quad.\]
 Now it is obvious that if $\xi$ is inside the unit disc, then $\ind\left(z-\xi\right)$
is positive (draw a picture!) --- this is not what we want. If $\xi$
is outside the disc, then $\ind\left(z-\xi\right)$ is harmonic (as
a function of $\xi$) and therefore has no maxima in any open subset
outside the disc, thus we need to investigate only the limiting behavior
as $\xi$ tends to infinity and to the unit circle.%
\footnote{The fact that the supremum occurs when $\xi$ converges to the unit
circle from outside can also be deduced by elementary geometric arguments.%
} Both can be calculated explicitly. When $\xi\rightarrow\infty$ clearly
$\ind\left(z-\xi\right)\rightarrow0.$ On the other hand, a simple
calculation shows that \begin{equation}
\lim_{\xi\rightarrow e^{2\pi is}}\ind\left(z-\xi\right)=\begin{cases}
\pi a-\pi & \textrm{if }s\in(0,a)\\
\pi a & \textrm{if }s\not\in[0,a]\end{cases}\label{eq:limind1}\end{equation}
 where the limit is taken from the outside of the unit circle. A geometric
note may be due here: the limit $\lim_{\xi\rightarrow e^{2\pi is}}\left(\widetilde{\arg}\left(e^{2\pi ia}-\xi\right)-\widetilde{\arg}\left(1-\xi\right)\right)$,
viewed as a function of $a$, for $s\in I$ is not continuous ---
it has a $-\pi$ jump discontinuity at $s$, this gives the $-\pi$
factor in (\ref{eq:limind1}).

A similar calculation shows that if $\xi\rightarrow1$ or $\xi\rightarrow e^{2\pi ia}$
then the limit of $\ind\left(z-\xi\right)$ depends on the angle of
approach, but in any case is bounded from below by $\pi a-\pi$. These
three cases ($\xi$ inside the disc, outside, and approaching the
boundary) totally give \begin{equation}
\ind\left(z-\xi\right)\geq\pi\left(a-1\right)\quad\forall\xi\in\mathbb{C}\setminus\left\{ e^{2\pi is}:s\in I\right\} \quad.\label{eq:ind1lt}\end{equation}
 The case $b=1$ is now immediate from (\ref{eq:indQsmall}) and (\ref{eq:ind1lt}).
The case $b>1$ is not much harder. Let $\tau=e^{2\pi i/b}$. If $\xi$
is a zero of $Q$ then so are $\tau\xi,\tau^{2}\xi,\dots,\tau^{b-1}\xi$
so we can rearrange $\xi_{j}$ such that $\xi_{j+K}=\tau\xi_{j}$.
In particular this gives the following representation using the first
$K$ zeros: \[
Q\left(z\right)=c\prod_{j=1}^{K}\left(z^{b}-\xi_{j}^{b}\right)\quad.\]
 Now the analog of (\ref{eq:ind1lt}) is\begin{equation}
\ind\left(z^{b}-\xi^{b}\right)\geq\pi\left(ba-1\right)\label{eq:indbgt}\end{equation}
 which holds for $a<1/b$ --- we will justify this assumption later.
The reasoning is similar: if $\xi$ is inside the unit disc so are
all the conjugates and the index is positive. For $\xi$ outside the
disc the index is harmonic in $\xi$ so we need to investigate only
$\xi$ converging to the circle, but as $a<1/b$ only one of the conjugates
of $\xi$ may converge to the arc $\left\{ e^{2\pi is}:s\in\left(0,a\right)\right\} $.
Summing $b$ copies of (\ref{eq:limind1}) gives (\ref{eq:indbgt}).
Summing the $K$ copies of (\ref{eq:indbgt}) gives (with (\ref{eq:indQsmall}))
\begin{equation}
K\pi\left(ba-1\right)\leq\ind Q\leq\pi-2\pi Na\label{eq:indfinal}\end{equation}
 which, again, holds for $a<1/b$. As according to our conditions
$2N/b>1$, for $a=\frac{1}{b}-\epsilon$ with $\epsilon>0$ sufficiently
small (\ref{eq:indfinal}) cannot hold, thus no such $f$ can have
an interval of length $\frac{1}{b}-\epsilon$ with no zeros. This
justifies our a priori assumption $a<1/b$, and we can use (\ref{eq:indfinal})
to get the largest $a$ allowed, proving (\ref{eq:MltM}). 

We turn to the proof of \begin{equation}
M\left(S\right)\geq\frac{K+1}{bK+2N}\quad.\label{eq:Mgt}\end{equation}
 The first step will be to construct an extremal function $f$ which
will only satisfy $f\geq0$ (i.e.~will have zeros in the interval
$I=\left[0,a\right]$ where $a=\frac{K+1}{bK+2N}$). The analysis
above suggests that we try $f=\real F$, $F\left(t\right)=e^{2\pi iNt}Q\left(e^{2\pi it}\right),$
for $Q$ satisfying (\ref{eq:Q_is_prod}), with $K$ zeros on the
arc $\left\{ e^{2\pi is}:s\in I\right\} $ (the other zeros are the
conjugates, which are of course also an the unit circle but not on
the arc). It turns out that it is enough to space the zeros evenly,
namely define \[
Q\left(z\right):=c\prod_{j=1}^{K}\left(z^{b}-\tau^{jb}\right)\]
 where $\tau:=e^{2\pi i\eta}$ with $\eta:=\frac{1}{bK+2N}$, and
$c$ is some constant such that $F\left(0\right)=Q\left(1\right)=-i$.
Now $\arg F$ can be calculated explicitly, where this time $\arg$
is in the usual sense, i.e.~a discontinuous function between $-\pi$
and $\pi$. Elementary geometry shows that \[
\arg\left(e^{i\varphi}-e^{i\nu}\right)\equiv\begin{cases}
\frac{\nu+\varphi-\pi}{2} & \textrm{if }\nu>\varphi\\
\frac{\nu+\varphi+\pi}{2} & \textrm{if }\nu<\varphi\end{cases}\mod2\pi\]
 Summing this for all the roots we get that on each segment $\left(k\eta,\left(k+1\right)\eta\right)$
$\arg F$ is a linear function with derivative $\pi\left(bK+2N\right)$.
Thus, on $\left(0,\eta\right)$ it is increasing from $-\pi/2$ to
$\pi/2$, at $\eta$ it has a jump of $-\pi$, then on $\left(\eta,2\eta\right)$
it is again linear with the same slope, so it again rises from $-\pi/2$
to $\pi/2$, then has another jump of $-\pi$ and so on. This argument
works till $\left(K+1\right)\eta=a$. Thus we get $\arg F\in\left[-\pi/2,\pi/2\right]$
on $I$ and therefore $f\geq0$ on $I$ having its only zeros at the
points $\left\{ 0,\eta,2\eta,\dots,K\eta,a\right\} $. 

Taking $f\left(t\right)+f\left(t+\epsilon\right)$ (whose Fourier
transform is also supported on $S$) for $\epsilon$ small enough
one gets a strictly positive function on any interval interior to
$I$. This proves (\ref{eq:Mgt}) and the theorem. 
\end{proof}

\section{Corollaries and remarks}

Let's compare the results of theorems 1 and 2. A simple calculation
shows that \[
D(S)=\sum_{n=N}^{N+K}\frac{1}{2n}=\frac{K+1}{2N+K}+O\left(\frac{K^{2}}{N^{3}}\right)\quad.\]
It is interesting to note that the true value, i.e.~$M(S)=\frac{K+1}{2N+K}$
is the unique approximation of $D(S)$ by a rational function of order
1 with this error term.

The following fact is a corollary of Theorem 2: 

\begin{cor*}
For every interval $I\subset[0,1]$ with $|I|<1$ there exists some
constant $\alpha$ such that for all $N\in\mathbb{N}$ there exists
a real trigonometric polynomial $f$ with \[
\spec f\subset[-\alpha N,-N]\cup[N,\alpha N]\]
such that $f>0$ on $I$.
\end{cor*}
\begin{proof}
Let $\delta:=1-|I|.$ Let $\alpha:=1+4/\delta.$ Apply Theorem 2 for
$N,$ $K:=\left\lfloor (\alpha-1)N\right\rfloor $ ($\left\lfloor \cdot\right\rfloor $
denoting as usual the integer value) and $\epsilon:=\delta/2.$ Then
as in this case \[
\left(\frac{K+1}{2N+K}-\epsilon\right)-|I|\geq\frac{1+\delta(N-\frac{1}{2})}{2N+K}>0,\]
 there exists a real trigonometric polynomial $f_{1}$ with $\spec\left(f_{1}\right)\subset\left[-\alpha N,-N\right]\cup\left[N,\alpha N\right]$
s.t.~$f_{1}>0$ on an interval of length bigger than $|I|$ (namely,
on the interval $\left[0,\frac{K+1}{2N+K}-\epsilon\right]$). If $I\subset\left[0,\frac{K+1}{2N+K}-\epsilon\right]$
we are done. Otherwise, we need another shift and $f\left(t\right):=f_{1}\left(t-a\right)$
will do, where $I=\left(a,b\right).$
\end{proof}
Here is a simple but slightly disappointing corollary to theorem \ref{thm:Mlinprog}:

\begin{cor*}
$M$ and $D$ may be significantly different even when they are both
small.
\end{cor*}
\begin{proof}
Take $b=K$ and $N=K^{2}$. Then easily\[
M(S)=\frac{1}{3}K^{-1}+O(K^{-2})\quad D(S)=(\log2)K^{-1}+O(K^{-3/2}).\qedhere\]

\end{proof}

\end{document}